\documentclass[11pt,reqno]{amsart}
\usepackage[usenames]{color}
\usepackage{amsmath,verbatim,graphicx,epstopdf,enumerate}
\newcommand{\D}{\mathrm{d}}

\newcommand{\lb}{\left(}

\newcommand{\rb}{\right)}
\newcommand{\PD}{\partial}

\newcommand{\Rb}{\mathbb{R}}

\newcommand{\Beq}{\begin{equation}}
	\newcommand{\Eeq}{\end{equation}}
\newcommand{\beq}{\begin{equation*}}
	\newcommand{\eeq}{\end{equation*}}
\newcommand{\bal}{\begin{align}}
	\newcommand{\eal}{\end{align}}

\newcommand{\bp}{\begin{prob}}
	\newcommand{\ep}{\end{prob}}
\newcommand{\bpr}{\begin{proof}}
	\newcommand{\epr}{\end{proof}}

\newcommand{\bel}[1]{\begin{equation}\label{#1}}
	\newcommand{\ee}{\end{equation}}

\newcommand{\rr}{\mathbb{R}}

\newtheorem{theorem}{Theorem}[section]

\newtheorem{lemma}[theorem]{Lemma}
\newtheorem{prop}[theorem]{Proposition}

\theoremstyle{definition}

\title[An Inverse problem with under-determined data]{Unique determination of the damping coefficient in the wave equation using point source and receiver data}
\author[Vashisth]{Manmohan Vashisth}
\address{Beijing Computational Science Research Center, Beijing 100193, China.
	\newline
	\indent E-mail: {\tt manmohanvashisth@gmail.com}}

\begin{document}
\maketitle

\begin{abstract}
In this article, we consider the inverse problems of determining the damping coefficient appearing in the wave equation. We prove the unique determination of the coefficient  from the data coming from a single coincident source-receiver pair. Since our problem is under-determined, so some extra assumption on the coefficient is required to prove the uniqueness.
\end{abstract}\vspace{2mm}

\vspace{.5mm}


\section{Introduction}
  We consider the following initial value problem (IVP),
\begin{equation}\label{wave equation with time derivative}
\begin{aligned}
(\Box- q(x)\partial_{t})u(x,t)&=\delta(x,t) \ \ \ \ (x,t)\in\rr^{3}\times\rr \\
u(x,t)|_{t<0}&=0\ \ \ \ \ \ \ \  \ \  \ \ \ \ \ \ \ \ \  x\in\rr^{3}
\end{aligned}
\end{equation}
where $\Box:=\PD_{t}^{2}-\Delta_{x}$ denotes the wave operator and the coefficient $q\in C^{\infty}(\Rb^{3})$ is known as damping coefficient. In this paper, we  study the problem of determination of coefficient $q$ appearing in \eqref{wave equation with time derivative} from the knowledge of solution measured at a single point for a certain period of time. We are interested in  the uniqueness of determination of coefficients $q$ from the knowledge of $u(0,t)$ for $t\in[0,T]$ with $T>0$ in Equation \eqref{wave equation with time derivative}. The problem studied here is motivated by geophysics, where geophysicists wish to determine the properties of earth structure by sending the waves from the surface of the earth and measuring the corresponding scattered responses (see \cite{Bube and Burridge 1-D,Symes} and references therein). Since the coefficient to be determined here depends on three variables while the given data depends on one variable as far as  the parameter count is concerned, the problem studied  here is under-determined. Thus some extra assumptions on coefficient $q$ are required in order to make the inverse problem solvable. We prove the uniqueness result for the radial coefficient.

There are several results related to the inverse problems for the wave equation with point source. We list them here. Romanov in \cite{Romanov} considered the problem for determining the  damping and potential coefficient in the wave equation with point source and proved unique determination of these coefficients by measuring the solution on a set containing infinite points. In \cite{Rakesh inversion of spherically symmetric potential} the problem of determining the radial potential from the knowledge of solution measured on a unit sphere for some time interval is studied. Rakesh and Sacks in \cite{Rakesh and Sacks uniqueness angular controlled potential}  established the  uniqueness for angular controlled potential in the wave equation from the knowledge of solution and its radial derivative measured on a unit sphere. In the above mentioned works the measurement set is an infinite set. Next we mention the work where uniqueness is established from the measurement of solution at a single point. Determination of the potential from the data coming from a single coincident source-receiver pair is considered in \cite{Rakesh single-coincident source-receiver pair} and the uniqueness result is established for the potentials which are either radial with respect a point different from source location or the potentials which are comparable. Recently author in \cite{Manmohan} extended the result of \cite{Rakesh single-coincident source-receiver pair} to a separated point source and receiver data. To  the best of our understanding,  very few results exist in the literature involving the recovery of the damping coefficient from point source and receiver data.  Our result, Theorem \ref{Damping Theorem 2}, is work in this direction. In the 1-dimensional inverse problems context, several results exist involving the uniqueness of recovery of the coefficient which depends on the space variable corresponding to  the first order derivative; see \cite{Ning_Yamamoto,Rakesh impedance inversion,Rakesh and Sacks impedance inversion,Rakesh Webster's equation,Romanov problem of determining the two coefficients,Sondhi}. We refer to \cite{Sacks_vocal_tract,Burridge integral equation,Levrentev Romanov and Shishatski book ill-posed,Rakesh layered medium point source,Rakesh and Uhlmann inverse backscattering} and references therein for more works related to the point source inverse problems for the wave equation.

We now state the main results of this article.

\begin{theorem}\label{Damping Theorem 2}
	Suppose $q_{i}(x)\in C^{\infty}(\rr^{3})$, $i=1,2$ with $q_{i}(x)=A_{i}(|x|)$ for some $C^{\infty}$ function $A_{i}$ on $[0,\infty)$. Let $u_{i}$ be the solution of the IVP
	\begin{align}\label{wave equation with time derivative2}
	\begin{aligned}
	(\Box- q_{i}(x)\partial_{t})u_{i}(x,t)&=\delta(x,t) \ \ \ \ (x,t)\in\rr^{3}\times\rr \\
	u_{i}(x,t)|_{t<0}&=0\ \ \ \ \ \ \ \  \ \  \ \quad \  x\in\rr^{3}.
	\end{aligned}
	\end{align}
	If $u_{1}(0,t)=u_{2}(0,t)$ for all $t\in[0,T]$ for some $T>0$, then $q_{1}(x)=q_{2}(x)$  for all $x$ with $|x|\leq T/2$, provided $q_{1}(0)=q_{2}(0)$.
\end{theorem}

The proof of the above theorem is based on an integral identity derived using the solution to an adjoint problem as used in \cite{Santosa Syemes} and \cite{Stefanov uniqueness}.  This idea was used in \cite{Eemeli,Rakesh and Uhlmann inverse backscattering,Manmohan} as well.

The article is organized as follows. In Section \ref{Prelim}, we state the existence and uniqueness results for the solution of  Equation \eqref{wave equation with time derivative},  the proof of which is given in \cite{Friedlander book,Levrentev Romanov and Shishatski book ill-posed,Romanov book integral geometry}.  Section \ref{Proof of Th 2} contains the proof of  Theorem \ref{Damping Theorem 2}.

\section{Preliminaries}\label{Prelim}
\begin{prop}\cite[pp.139,140]{Friedlander book}\label{Damping Prop}
	Suppose $q\in C^{\infty}(\Rb^{3})$
	and $u(x,t)$ satisfies the following initial value problem
	\begin{equation}\label{Hyperbolic PDE}
	\begin{aligned}
	Pu(x,t):=(\Box -q(x)\partial_{t})u(x,t)&= \delta(x,t),  \ (x,t)\in\rr^{3}\times\rr\\
	u(x,t)|_{t<0}\ &=\ 0, \quad \quad \quad \ \ \  x\in \rr^{3}
	\end{aligned}
	\end{equation}
	then $u(x,t)$ is given by 
	\begin{equation}\label{Fundamental solution to hyperbolic equation}
	u(x,t) =\frac{R(x,t)\delta(t-|x|)}{4\pi|x|} + v(x,t)
	\end{equation}
	where $v(x,t)=0$ for $t<|x|$ and in the region $t>|x|$, $v(x,t)$ is a $C^{\infty}$
	solution of the characteristic boundary value problem (Goursat Problem)
	\begin{equation}\label{Goursat problem}
	\begin{aligned}
	Pv(x,t)&= 0, \ for \ t>|x|\\
	v(x,|x|)&=-\frac{R(x,|x|)}{8\pi} \int\limits_{0}^{1}\frac{PR(sx,s|x|)}{R(sx,s|x|)} ds ,\ \ \  \forall x\in\rr^{3}
	\end{aligned}
	\end{equation}
and 	$R(x,t)$ is given by \cite[pp. 134]{Friedlander book}
	\begin{equation}\label{Expression for r}
	R(x,t)=\exp\left(-\frac{1}{2}\int\limits_{0}^{1}q(sx)t ds\right).
	\end{equation}
	
\end{prop}

\renewcommand{\thesection}{\large 3}
\section{Proof of Theorem \ref{Damping Theorem 2}}\label{Proof of Th 2}
In this section, we prove Theorem \ref{Damping Theorem 2}.   We will first prove an integral identity which will be used to prove our main result.
\begin{lemma}\label{Damping integral identity}
	Let $u_{i}(x,t)$ for $i=1,2$ be the solution to Equation \eqref{wave equation with time derivative2}. Then the following integral identity holds for all $\sigma\geq 0$
	\begin{equation}\label{Integral Identity for time derivative case}
	\int\limits_{\rr^{3}}\int\limits_{\rr}q(x)\partial_{t}u_{2}(x,t)u_{1}(x,2\sigma-t)dt dx =u(0,2\sigma)
	\end{equation}
	where $q(x):=q_{1}(x)-q_{2}(x)$ 
	and $u(x,t)=(u_{1}-u_{2})(x,t)$.\\
	\begin{proof} 
		Here we have  $u$ satisfies the following IVP 
		\begin{equation}\label{equation for u in time derivative case}
		\begin{aligned}
		\Box u(x,t)-q_{1}(x)\partial_{t}u(x,t)&=q(x)\partial_{t}u_{2}(x,t)\  \ \ \ \ (x,t)\in\rr^{3}\times\rr\\
		u(x,t)|_{t<0}&=0\ \ \ \ \ \ \ \ \ \ \  \ \ \ \ \ \ \ \ \ \ \ \ \ \ \  x\in\rr^{3}.
		\end{aligned}
		\end{equation}
		Multiplying Equation \eqref{equation for u in time derivative case} by $u_{1}(x,2\sigma-t)$ and integrating over $\rr^{3}\times\rr$, we have 
		\begin{align*}
		\begin{aligned}
	\int\limits_{\rr^{3}}\int\limits_{\rr}q(x)\partial_{t}u_{2}(x,t)u_{1}(x,2\sigma-t)dtdx &=	\int\limits_{\rr^{3}}\int\limits_{\rr}\left(\Box u(x,t) -q_{1}(x)\partial_{t}u(x,t)\right)u_{1}(x,2\sigma-t) dt dx\\
	&=\int\limits_{\Rb^{3}}\int\limits_{\Rb}u(x,t)\lb\Box u(x,t) -q_{1}(x)\partial_{t}u_{1}(x,2\sigma-t)\rb\D x \D t
		\end{aligned}
		\end{align*}
		where in the last step above we have used integration by parts and the properties of $v$ in Proposition \ref{Damping Prop}. Thus finally using the fact that $u_{1}$ is solution to \eqref{wave equation with time derivative2}, we get 
		\begin{align*}
		\begin{aligned}
		\int\limits_{\rr^{3}}\int\limits_{\rr}q(x)\partial_{t}u_{2}(x,t)u_{1}(x,2\sigma-t)dtdx=u(0,2\sigma);\ \mbox{for all $\sigma\geq 0$}.
		\end{aligned}
		\end{align*}
		 This completes the proof of the lemma.
	\end{proof}
\end{lemma}

Using Lemma \ref{Damping integral identity}   and  the fact that $u(0,t)=0$ for all $t\in[0,T]$, we see that
\begin{align*}
\int\limits_{\rr^{3}}\int\limits_{\rr}q(x)\partial_{t}u_{2}(x,t)u_{1}(x,2\sigma-t)dtdx=0; \mbox{ for all }\sigma\in[0,T/2].
\end{align*}
Now using Equation $\eqref{Fundamental solution to hyperbolic equation}$, we get 
\begin{align*}
\begin{aligned}
&\int\limits_{\rr^{3}}\int\limits_{\rr}q(x)\partial_{t}\Bigg(\frac{R_{2}(x,t)\delta(t-|x|)}{4\pi |x|}+v_{2}(x,t)\Bigg) \Bigg( \frac{R_{1}(x,2\sigma-t)\delta(2\sigma-t-|x|)}{4\pi |x|}\Bigg)dtdx\\
&+ \int\limits_{\rr^{3}}\int\limits_{\rr}q(x)v_{1}(x,2\sigma-t)\partial_{t}\Bigg(\frac{R_{2}(x,t)\delta(t-|x|)}{4\pi |x|}+v_{2}(x,t)\Bigg)dtdx=0.
\end{aligned}
\end{align*}
This gives
\begin{equation}\label{Explicit form of integral identity in time derivative case}
\begin{aligned}
&\underbrace{\int\limits_{\rr^{3}}\int\limits_{\rr} \frac{q(x)\partial_{t}R_{2}(x,t)R_{1}(x,2\sigma-t)\delta(t-|x|)\delta(2\sigma-t-|x|)}{16\pi^{2}|x|^{2}}dtdx}_{I_{1}}\\ &+\underbrace{\int\limits_{\rr^{3}}\int\limits_{\rr}\frac{q(x)R_{2}(x,t)R_{1}(x,2\sigma-t)\partial_{t}\delta(t-|x|)\delta(2\sigma-t-|x|)}{16\pi^{2}|x|^{2}}dtdx}_{I_{2}}\\
&+\underbrace{\int\limits_{\rr^{3}}\int\limits_{\rr}q(x)\partial_{t}\Big(\frac{R_{2}(x,t)\delta(t-|x|)}{4\pi |x|}\Big)v_{1}(x,2\sigma-t)dtdx}_{I_{3}}\\
&+\underbrace{\int\limits_{\rr^{3}}\int\limits_{\rr}\frac{q(x)\partial_{t}v_{2}(x,t)R_{1}(x,2\sigma-t)\delta(2\sigma-t-|x|)}{4\pi |x|}dtdx}_{I_{4}}\\
&+\underbrace{\int\limits_{\rr^{3}}\int\limits_{\rr}q(x)\partial_{t}v_{2}(x,t)v_{1}(x,2\sigma-t)dtdx}_{I_{5}}=0;\mbox{ for all }  \sigma\in [0,T/2].
\end{aligned}
\end{equation}
In a compact form, this can be written as
\begin{equation}\label{integral identity in terms of Ij}
I_{1}+I_{2}+I_{3}+I_{4}+I_{5}=0.
\end{equation}
Next we simplify each $I_{j}$ with $j=1,2,.....,5$. We will use the fact that $v_{i}(x,t)=0$ for $t<|x|$. 

We have 
\begin{align*}
I_{1}&=\int\limits_{\rr^{3}}\int\limits_{\rr} \frac{q(x)\partial_{t}R_{2}(x,t)R_{1}(x,2\sigma-t)\delta(t-|x|)\delta(2\sigma-t-|x|)}{16\pi^{2}|x|^{2}}dt dx\\
&=\int\limits_{|x|=\sigma}\frac{q(x)\partial_{t}R_{2}(x,|x|)R_{1}(x,|x|)}{16\pi^{2}|x|^{2}}dS_{x}\\
&=-\int\limits_{|x|=\sigma}\frac{q(x)R_{1}(x,|x|)R_{2}(x,|x|)}{32\pi^{2}|x|^{2}}\Big(\int\limits_{0}^{1}q_{2}(sx)ds\Big)dS_{x}.
\end{align*}
Next we simplify the integral $I_{2}$. We use the following formula \cite[Page 231, Eq.(10)]{Gelfand}
\begin{equation}\label{Dirac derivative}
\int \delta'(r-|x|)\varphi dx=\frac{-1}{|x|^{2}}\int\limits_{|x|=r} \frac{\partial }{\partial r} \left (\varphi r^{2}\right) dS_{x}.
\end{equation}
Note that from  this formula, by a change of variable, we have
\begin{equation}\label{Dirac derivative 2}
\int \delta'(2r-2|x|)\varphi dx=\frac{-1}{2|x|^{2}}\int\limits_{|x|=r} \frac{\partial }{\partial r} \left (\varphi r^{2}\right) dS_{x}.
\end{equation}
Now
\begin{align*}
I_{2}&=\int\limits_{\rr^{3}}\int\limits_{\rr}\frac{q(x)R_{2}(x,t)R_{1}(x,2\sigma-t)\partial_{t}\delta(t-|x|)\delta(2\sigma-t-|x|)}{16\pi^{2}|x|^{2}}dtdx\\
&=\int\limits_{\rr^{3}}\int\limits_{\rr}\frac{q(x)R_{2}(x,t)R_{1}(x,2\sigma-t)\delta^{'}(t-|x|)\delta(2\sigma-t-|x|)}{16\pi^{2}|x|^{2}}dtdx\\
&=\int\limits_{\rr^{3}}\frac{q(x)R_{2}(x,2\sigma-|x|)R_{1}(x,|x|)\delta'(2\sigma-2|x|)}{16\pi^{2}|x|^{2}}dx\\
&=-\frac{1}{32\pi^{2}\sigma^{2}}\int\limits_{|x|=\sigma}\frac{\partial}{\partial{r}}\{q(x)R_{1}(x,|x|)R_{2}(x,2\sigma-|x|)\}dS_{x}.
\end{align*}
In the last step above, we used Equation \eqref{Dirac derivative 2}.

Next we have 
\begin{align*}
I_{3}&=\int\limits_{\rr^{3}}\int\limits_{\rr}q(x)\partial_{t}\Big(\frac{R_{2}(x,t)\delta(t-|x|)}{4\pi |x|}\Big)v_{1}(x,2\sigma-t)dxdt.\\
\intertext{{We can view the derivative above as a limit of the difference quotients in the distribution topolgy \cite[pp.48]{Friedlander-Distributions-book}. Combining this with the fact that $v_{1}$ is $C^{2}$ in $\{(x,t):|x|\leq t\}$, we get,}}
I_{3}&=-\int\limits_{\rr^{3}}\int\limits_{\rr}q(x)\frac{R_{2}(x,t)\delta(t-|x|)}{4\pi |x|}\partial_{t}\Big(v_{1}(x,2\sigma-t)\Big)dxdt\\
&=\int\limits_{\rr^{3}}\frac{q(x)R_{2}(x,|x|)\partial_{t}v_{1}(x,2\sigma-|x|)}{4\pi |x|}dx.\\
\intertext{Again using the fact that $v_{1}(x,t)=0$ for $t<|x|$, we get,}
&I_{3}=\int\limits_{|x|\leq \sigma}\frac{q(x)R_{2}(x,|x|)\partial_{t}v_{1}(x,2\sigma-|x|)}{4\pi |x|}dx.
\end{align*}
Next we simplify $I_{4}$. Similiar to $I_{3}$, we have
\begin{align*}
I_{4}&=\int\limits_{\rr^{3}}\int\limits_{\rr}\frac{q(x)\partial_{t}v_{2}(x,t)R_{1}(x,2\sigma-t)\delta(2\sigma-t-|x|)}{4\pi |x|}dtdx\\
&=\int\limits_{|x|\leq \sigma}\frac{q(x)R_{1}(x,|x|)\partial_{t}v_{2}(x,2\sigma-|x|)}{4\pi|x|}dx.
\end{align*}
Finally, we have
\begin{align*}
I_{5}&=\int\limits_{\rr^{3}}\int\limits_{\rr}q(x)\partial_{t}v_{2}(x,t)v_{1}(x,2\sigma-t)dtdx\\
&=\int\limits_{|x|\leq\sigma}\int\limits_{|x|}^{2\sigma-|x|}q(x)\partial_{t}v_{2}(x,t)v_{1}(x,2\sigma-t)dtdx.
\end{align*}
Now, we use the fact that $q_i$ is a radial function, that is, $q_{i}(x)=A_{i}(|x|)$.  Then note that
\[
R_{i}(x,|x|)= \exp\left(-\frac{|x|}{2} \int\limits_{0}^{1} q_{i}(sx)ds\right)=\exp\left(-\frac{|x|}{2} \int\limits_{0}^{1} A_{i}(s|x|)ds\right)
\]
is also radial. For simplicity, we denote $R(x,|x|)$ by $R(|x|)$.

With this, we have
\[
I_{1}=-\frac{A(\sigma)R_{1}(\sigma)R_{2}(\sigma)}{8\pi}\int\limits_{0}^{1} A_{2}(s\sigma) ds.
\]
Next we consider $I_{2}$. First let us consider the derivative:
\[
D_{r}:=\frac{\partial}{\partial r}\left( A(r) R_{1}(x,r) R_{2}(x,2\sigma-r)\right).
\]
After a routine calculation, we get,
\begin{align*}
D_{r}&= A'(r) R_{1}(x,r)R_{2}(x,r) -\frac{1}{2} A(r)^{2} R_{1}(x,r)R_{2}(x,2\sigma-r) \\
&-\sigma A(r)R_{1}(x,r)R_{2}(x,2\sigma-r) \int\limits_{0}^{1} A_{2}'(rs)s ds\\
&= A'(r) R_{1}(x,r)R_{2}(x,r) -\frac{1}{2} A(r)^{2} R_{1}(x,r)R_{2}(x,2\sigma-r) \\
&-A(r)R_{1}(x,r)R_{2}(x,2\sigma-r)\left[\frac{\sigma}{r}\left(A_{2}(r)-\int\limits_{0}^{1}A_{2}(rs) ds\right)\right].
\end{align*}
On $|x|=\sigma$, we have
\begin{align*}
D_{r}|_{|x|=\sigma}&=R_{1}(\sigma)R_{2}(\sigma)\left[A'(\sigma)-\frac{1}{2} A(\sigma)^{2}-A(\sigma)A_{2}(\sigma)+A(\sigma)\int\limits_{0}^{1}A_{2}(s\sigma) ds\right]\\
&=R_{1}(\sigma)R_{2}(\sigma)\left[A'(\sigma) -\frac{1}{2}A(\sigma)(A_{1}+A_{2})(\sigma)+A(\sigma)\int\limits_{0}^{1}A_{2}(s\sigma) ds\right].
\end{align*}
Hence
\begin{align*}
I_{2}&=-\frac{1}{8\pi}\left(R_{1}(\sigma)R_{2}(\sigma)\left[A'(\sigma) -\frac{1}{2}A(\sigma)(A_{1}+A_{2})(\sigma)+A(\sigma)\int\limits_{0}^{1}A_{2}(s\sigma) ds\right]\right).
\end{align*}
Let us denote 
\[
\tilde{A}(\sigma)=A(\sigma)R_{1}(\sigma)R_{2}(\sigma).
\]
Then
\begin{align*}
I_{2}=-\frac{1}{8\pi} \frac{d}{d\sigma}\tilde{A}(\sigma)-\frac{1}{8\pi}\tilde{A}(\sigma)\int\limits_{0}^{1}A_{2}(s\sigma)ds.
\end{align*}
Therefore
\[
I_{1}+I_{2}=-\frac{1}{8\pi} \left(2\tilde{A}(\sigma)\int\limits_{0}^{1}A_{2}(s\sigma) ds+\frac{d}{d\sigma}\tilde{A}(\sigma)\right).
\]
Considering the following integrating factor for $I_{1}+I_{2}$ 
\[
\exp\left(2\int\limits_{0}^{\sigma}\int\limits_{0}^{1}A_{2}(ts)dt ds\right),
\]
we have
\[
I_{1}+I_{2}=-\frac{1}{8\pi}\exp\left(-2\int\limits_{0}^{\sigma}\int\limits_{0}^{1}A_{2}(ts)dt ds\right)\frac{d}{d\sigma}\left[\exp\left(2\int\limits_{0}^{\sigma}\int\limits_{0}^{1}A_{2}(ts)dt ds\right)\tilde{A}(\sigma)\right].
\]

Now from Equation \eqref{integral identity in terms of Ij},     we have
\begin{equation}\label{integral equation in radial case}
\begin{aligned}
& \frac{1}{8\pi}\frac{d}{d\sigma}\left[\tilde{A}(\sigma)\exp\left(2\int\limits_{0}^{\sigma}\int\limits_{0}^{1}A_{2}(st)dsdt\right)\right]\\
&\ \ =\exp\left(2\int\limits_{0}^{\sigma}\int\limits_{0}^{1}A_{2}(st)dsdt\right)\Bigg[\int\limits_{|x|\leq \sigma}\frac{q(x)R_{2}(x,|x|)\partial_{t}\{R_{1}v_{1}\}(x,2\sigma-|x|)}{4\pi |x|}dx\\
&\ \ \ \ +\int\limits_{|x|\leq\sigma}\frac{q(x)R_{1}(x,|x|)\partial_{t}v_{2}(x,2\sigma-|x|)}{4\pi|x|}dx\\
&\ \ \ \ +\int\limits_{|x|\leq\sigma}\int\limits_{|x|}^{2\sigma-|x|}q(x)\partial_{t}v_{2}(x,t)v_{1}(x,2\sigma-t)dtdx\Bigg]\mbox{ for all } \sigma\in[0,T/2].
\end{aligned}
\end{equation}
Integrating on both sides with respect to $\sigma$ under the assumption  that $\tilde{A}(0)=0$, we get 
\begin{align*}
&\exp\left(\int\limits_{0}^{\tilde{\sigma}}\int\limits_{0}^{1}2 A_{2}(st)dsdt\right)\tilde{A}(\tilde{\sigma})\\
&=\int\limits_{0}^{\tilde{\sigma}}\exp\left(\int\limits_{0}^{\sigma}\int\limits_{0}^{1}2 A_{2}(st)dsdt\right)\Bigg\{\int\limits_{|x|\leq \sigma}\frac{q(x)R_{2}(x,|x|)\partial_{t}v_{1}(x,2\sigma-|x|)}{4\pi |x|}dx\\
&+\int\limits_{|x|\leq\sigma}\frac{q(x)R_{1}(x,|x|)\partial_{t}v_{2}(x,2\sigma-|x|)}{4\pi|x|}dx\\
&+\int\limits_{|x|\leq\sigma}\int\limits_{|x|}^{2\sigma-|x|}q(x)\partial_{t}v_{2}(x,t)v_{1}(x,2\sigma-t)dtdx\Bigg\} d\sigma, \  \mbox{ for all } \tilde{\sigma}\in[0,T/2].
\end{align*}
Now using the fact that $R_{i}'s$ are continuous, non-zero functions, and $v_{i}'s$ are continuous, we have the following inequality:
\begin{align*}
|\tilde{A}(\tilde{\sigma})|\leq C\int\limits_{0}^{\tilde{\sigma}}|\tilde{A}(r)| dr\mbox{ for all }\tilde{\sigma}\in[0,T/2].
\end{align*}
Now by Gronwall's inequality, we have $\tilde{A}(\sigma)=0$ for all $\tilde{\sigma}\in [0,T/2]$, which gives us $q_{1}(x)=q_{2}(x)$ for all $x\in\rr^{3}$ such that $|x|\leq T/2$. This completes the proof.

\section*{Acknowledgement}
The author would like to thank Dr. Venky Krishnan for useful discussions. He is supported by NSAF grant (No. U1530401).

\end{document}